\documentclass[12pt]{article}
\usepackage{bbm}
\usepackage{latexsym}
\usepackage{mathrsfs}
\usepackage{amsmath}
\usepackage{graphicx}
\usepackage{amssymb}
\usepackage{color}
\usepackage{amsthm}
\usepackage{cite}
\usepackage{indentfirst}
\usepackage{anysize}\marginsize{45mm}{45mm}{40mm}{50mm}
\usepackage{caption}
\usepackage{enumerate}

\renewcommand\thesection{\arabic{section}.\kern -.5em}

\newtheoremstyle{lemma}{\topsep}{\topsep}%
     {}
     {}
     {\bfseries}
     {}
     {0.1em}
     {\thmname{#1}\thmnumber{ #2}\thmnote{ #3}}
\theoremstyle{lemma}  

\newtheorem{theorem}{Theorem}              
\newtheorem{lemma}[theorem]{Lemma}

\numberwithin{equation}{section}

\usepackage{latexsym, bm}
\begin{document}
\title{Paired many-to-many 2-disjoint path cover of Johnson graphs\thanks{This research was partially supported by the National Natural Science Foundation of China (Nos. 11801061 and 62371094)}}

\author{Jinhao Liu and Huazhong L\"{u}\thanks{Corresponding author.}\\
{\small School of Mathematical Sciences, } \\
{\small University of Electronic Science and Technology of China,} \\
{\small Chengdu, Sichuan 610054, P.R. China}\\
{\small E-mail: 202321110137@std.uestc.edu.cn; lvhz@uestc.edu.cn}\\}
\date{}

\maketitle

\begin{abstract}
Given two 2 disjoint vertex-sets $S=\{u,x\}$ and $T=\{v,y\}$, a paired many-to-many 2-disjoint path cover joining S and T, is a set of two vertex-disjoint paths with endpoints $u,v$ and $x,y$, respectively, that cover every vertex of the graph. If the graph has a many-to-many 2-disjoint path cover for any two disjoint vertex-sets $S$ and $T$, then it is called paired 2-coverable. It is known that if a graph is paired 2-coverable, then it must be Hamilton-connected, but the reverse is not true. It has been proved that Johnson graphs $J(n,k)$, $0\le k\le n$, are Hamilton-connected by Brian Alspach in [Ars Math. Contemp. 6 (2013) 21--23]. In this paper, we prove that Johnson graphs are paired 2-coverable. Moreover, we obtain that another family of graphs $QJ(n,k)$ constructed from Johnson graphs by Alspach are also paired 2-coverable.

\vskip 0.1 in


\end{abstract}

\section{Introduction}
Let $[n]=\{1,2,\cdots,n\}$. The Johnson graph $J(n,k)$, $0\le k \le n$, is defined by letting vertices correspond to $k$-subsets of $[n]$. Two vertices are adjacent if their corresponding $k$-subsets have $k-1$ common elements. For simplicity, we denote $n \in u$ if the element $n$ belongs to the $k$-subset corresponding to the vertex $u$. The graphs $QJ(n,A)$ are defined as follows. Let $A=\{a_1,a_2,\dotsi ,a_m\}$ be a non-empty subset of $[n]$ such that the elements are listed in the order $a_1<a_2<\dotsi<a_m$. For each $a_i \in A$, we first take a copy of the Johnson graph $J(n,a_i)$. For each $i$, we add an edge between the vertex $u$ in $J(n,a_i)$ and the vertex $v$ in $J(n,a_{i+1})$ if $u$ is a subset of $v$. For simplicity, we say that a vertex of $J(n,a_i)$ in $QJ(n,A)$ lies in level $i$. A graph is {\em Hamilton-connected} if for any pair of distinct vertices $u,v$ there is a Hamilton path whose terminal vertices are $u$ and $v$. Let $G$ be a graph and we denote the vertex set of $G$ (resp. edge set of $G$) by $V(G)$ (resp. $E(G)$). We use $\langle u_1,u_2,\cdots, u_n\rangle$ to denote a path from $u_1$ to $u_n$ via $u_2,u_3,\cdots,u_{n-1}$ in order. A graph has a {\em paired many-to-many 2-disjoint path cover}, if given two disjoint vertex-sets $S=\{u,x\}$ and $T=\{v,y\}$, there are two vertex-disjoint paths with endpoints $u,v$ and $x,y$, respectively, that cover every vertex of the graph. If the graph has a many-to-many 2-disjoint path cover for any two disjoint vertex-sets $S$ and $T$ of the graph, it is called {\em paired 2-coverable}. For convenience, we denote the two paths by $P2C(u,v;x,y)$. In addition, $P2C(u,v;x,y)$ will be abbreviated as $P2C$ paths if the context is clear. And we always use $u,v,x$ and $y$ as endpoints of $P2C$ paths unless state otherwise.

\begin{figure}
	\centering
	\includegraphics[width=0.4\textwidth]{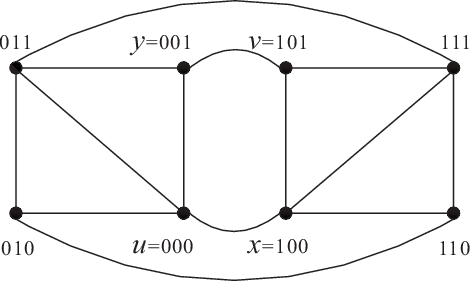}
     \caption{A Hamilton-connected graph which is not paired 2-coverable. }
     \label{ex}
\end{figure}

By choosing $u,v$ and $x,y$ appropriately in a graph $G$, i.e. $vy\in E(G)$, it can be easily to yield a Hamilton path from $u$ to $x$ by adding the edge $vy$ to $P2C(u,v;x,y)$. This implies that a graph which is paired 2-coverable must be Hamilton-connected. However, not all Hamilton-connected graphs are paired 2-coverable. Here is an example, as shown in Fig.~\ref{ex}. It is a $3$-dimensional hypercube by adding two edges $\{000,011\}$ and $\{100,111\}$, denoted by $G$. It is trivial to verify that $G$ is Hamilton-connected. Let $u=000$, $v=101$, $x=100$ and $y=001$. If there exists two paths $P$ and $Q$ of $P2C(u,v;x,y)$, the path $P$ from $u=000$ to $v=101$ must contain one of the edges $\{011, 111\}$ and $\{110, 111\}$. This implies that the neighbors of $x$ or $y$ are all contained in the path $P$ from $u$ to $v$. As a result, there is no path $Q$ with endpoints $x$ and $y$ disjoint from $P$.

Moreover, many-to-many 2-disjoint path cover for some well-known graphs have been studied~\cite{CHF,HZ,JI,JH,KJ}. Recently, Brian Alspach \cite{Alspach} showed that $J(n,k)$ and $QJ(n,A)$ are Hamilton-connected. We are interested in considering whether $J(n,k)$ and $QJ(n,A)$ are paired 2-coverable. In this paper, we prove that $J(n,k)$ is paired 2-coverable whenever $n \ge 4,n>k\ge 1$ by double induction. Finally, we prove that $QJ(n,A)$ is paired 2-coverable for all $n\ge 4$.

\vskip 0.0 in

\section{Main results}
To begin with, we present the following lemmas.

\begin{lemma}
    \cite{Alspach} $J(n,k)$ is Hamilton-connected for all $n\ge 1$.
\end{lemma}
\begin{lemma}
     \cite{Alspach} $QJ(n,A)$ is Hamilton-connected for all $n\ge 3$.
\end{lemma}

\begin{theorem}
    ~The complete graph $K_n$ is paired 2-coverable for all $n \ge 4$.
\end{theorem}
\noindent Proof. For any two pairs of distinct vertices $u,v$ and $x,y$, we let $\{a_1,a_2,\cdots ,a_{n-4}\}$ be the set of vertices of $K_n-\{u,v,x,y\}$. Thus, $\langle u,v\rangle$ and $\langle x,a_1,\cdots ,a_{n-4},y\rangle$ are $P2C$ paths of $K_n$. This completes the proof.$\qed$

\begin{theorem}
    ~ $J(n,k)$ is paired 2-coverable when $n \ge 4$ and $n>k\ge 1$.
\end{theorem}
\noindent Proof. It is obvious that  $J(n,1)$, $n\ge 2$, is isomorphic to the complete graph $K_n$ and that $J(n,k)$ and $J(n,n-k)$ are isomorphic via mapping a $k$-subset to its complement~\cite{CG}.
In the following, we prove this theorem by double induction. First we need to verify that $J(n,1)$ for $n \ge 4$, $J(4,2)$ and $J(k+1,k)$ for $k \ge 3$, are paired 2-coverable.
 Since $J(n,1)$, $n\ge 4$, is isomorphic to $K_n$ and $J(k+1,k)$, $k\ge 3$, is isomorphic to $K_{k+1}$, by Theorem 3, they are paired 2-coverable. For $J(4,2)$, each vertex corresponds to a $2$-subset of $\{1,2,3,4\}$. Since $J(4,2)$ is isomorphic to $K_{2,2,2}$, which is vertex and edge-transitive, we list all essentially distinct $P2C$ paths of $J(4,2)$ in the following table (Table 1).
 \begin{table}[h]\scriptsize
\centering
\begin{tabular}{ccccll} \hline
$u$ & $v$ & $x$ & $y$ & a path from $u$ to $v$ & a path from $x$ to $y$ \\
\hline
$\{1,2\}$ & $\{1,3\}$ & $\{2,3\}$ & $\{2,4\}$ &  $\langle\{1,2\},\{1,4\},\{1,3\}\rangle$ & $\langle\{2,3\},\{3,4\},\{2,4\}\rangle$ \\
$\{1,2\}$ & $\{2,4\}$ & $\{1,3\}$ & $\{2,3\}$ & $\langle\{1,2\},\{1,4\},\{2,4\}\rangle$ & $\langle\{1,3\},\{3,4\},\{2,3\}\rangle$ \\
 $\{1,2\}$ & $\{2,3\}$& $\{1,3\}$ & $\{2,4\}$ & $\langle\{1,2\},\{2,3\}\rangle$ & $\langle
\{1,3\},\{1,4\},\{3,4\},\{2,4\}\rangle$ \\
$\{1,2\}$ & $\{1,3\}$ & $\{2,4\}$ & $\{3,4\}$ &  $\langle\{1,2\},\{1,4\}\{1,3\}\rangle$ & $\langle\{2,4\},\{2,3\},\{3,4\}\rangle$ \\
$\{1,2\}$ & $\{2,4\}$ & $\{1,3\}$ & $\{3,4\}$ & $\langle\{1,2\},\{2,3\},\{2,4\}\rangle$ & $\langle\{1,3\},\{1,4\},\{3,4\}\rangle$ \\
 $\{1,2\}$ & $\{3,4\}$& $\{1,3\}$ & $\{2,4\}$ & $\langle\{1,2\},\{2,3\},\{3,4\}\rangle$ & $\langle
\{1,3\},\{1,4\},\{2,4\}\rangle$ \\
\hline
\end{tabular}
\caption{All essentially distinct $P2C$ paths of $J(4,2)$.}\label{t1}
\end{table}

   When considering $J(n,k)$, the induction hypotheses are: $J(m,k')$ is paired 2-coverable whenever $k' < k$ and $k' < m \le n$ or $J(m,k)$ is paired 2-coverable whenever $k < m <n$.

   If $k<n<2k$, then $n-k<k$ so that $J(n,n-k)$ is paired 2-coverable by induction. It follows that $J(n,k)$ is paired 2-coverable because $J(n,k)$ is isomorphic to $J(n,n-k)$.

   If $n \ge 2k$, let $u,v,x,y$ be four distinct vertices in $J(n,k)$. Next we prove that there are two paths of $P2C(u,v;x,y)$ in $J(n,k)$. Let $X$ be the induced subgraph of all vertices that do not contain the element $n$ in $J(n,k)$ and let $Y$ be the induced subgraph of all vertices that contain the element $n$ in $J(n,k)$. It is clear that $X$ is isomorphic to $J(n-1,k)$ and $Y$ is isomorphic to $J(n-1,k-1)$. We distinguish the following cases according to the number of endpoints that contain $n$.\\
   \textbf{Case 1}. All the four endpoints contain $n$. Thus, $u,v,x,y$ belong to $Y$, implying that there are two paths $P$ and $Q$ of $P2C(u,v;x,y)$ in $Y$ by induction. Let $ab$ be an edge on $P$. It is clear that $n-1>k$. Thus, replacing $n$ by an element, say $i(\neq n)$, that neither $a$ nor $b$ contains, we can choose $a'$ and $b'$ containing $i$ as neighbors of $a$ and $b$, respectively. Thus, there is a Hamilton path $R$ from $a'$ to $b'$ in $X$ by Lemma 1. Deleting the edge $ab$, adding the edges $aa',bb'$ and concatenating the path $R$. Hence, we obtain two paths of $P2C(u,v;x,y)$ in $J(n,k)$.\\
   \textbf{Case 2}. None of four endpoints contain $n$. Thus, $u,v,x,y$ belong to $X$ implying that there are two path $P$ and $Q$ of $P2C(u,v;x,y)$ in $X$ by induction. Let $ab$ be an edge on $P$. By using the approach analogous to Case 1, we can find a vertex $a'$ adjacent to $a$ in $Y$ and a vertex $b'$ adjacent to $b$ in $Y$. Thus, there is a Hamilton path $R$ from $a'$ to $b'$ in $Y$ by Lemma 1. Deleting the edge $ab$, adding the edges $aa',bb'$ and concatenating the path $R$. Hence, we obtain two paths of $P2C(u,v;x,y)$ in $J(n,k)$.\\
   \textbf{Case 3}. Exactly one endpoint contains $n$. Without loss of generality, we assume that only $u$ contains $n$. Choose a vertex $a$ that do not contain the element $n$ and different from $v,x$ or $y$. Thus, there are two paths of $P2C(a,v;x,y)$ in $X$ by induction. Let $b\ne u$ be a vertex adjacent to $a$ in Y. Thus, there is a Hamilton path $R$ from $u$ to $b$ in $Y$ by Lemma 1. Then adding the edge $ab$ and concatenating the path $R$, we obtain two paths of $P2C(u,v;x,y)$ in $J(n,k)$.\\
   \textbf{Case 4}. Exactly one endpoint does not contain $n$. Without loss of generality, we assume only $u$ does not contain $n$. Choose a vertex $a$ that contains the element $n$ and different from $v,x$ or $y$. Thus, there are two paths of $P2C(a,v;x,y)$ in $Y$ by induction. Let $b(\neq u)$ be a vertex in $X$ which is adjacent to $a$. Thus, there is a Hamilton path $R$ from $u$ to $b$ in $X$ by Lemma 1. Then adding the edge $ab$ and concatenating the path $R$, we obtain the two paths of $P2C(u,v;x,y)$ in $J(n,k)$.\\
   \textbf{Case 5}. Exactly two endpoints contain $n$. We may assume that each element of $[n]$ appears in exactly two endpoints. Otherwise, we can choose such an element (not $n$) in $[n]$ to replace $n$. Since each element of $[n]$ appears twice, and the four endpoints consist of $4k$ elements, $2n=4k$. That is $n=2k$.\\
   \textbf{Case 5.1}. $n$ is contained in $u$ and $v$, i.e. $u,v$ belong to $Y$ and $x,y$ belong to $X$. It follows from Lemma 1 that there is a Hamilton path $P$ from $u$ to $v$ in $Y$ (resp. $Q$ from $x$ to $y$ in $X$). Thus $P$ and $Q$ are the two paths of $P2C(u,v;x,y)$ in $J(n,k)$.\\
   \textbf{Case 5.2}. $n$ is contained in $v$ and $y$. We may assume that $u =\{1,2,\dotsi,k\},v = \{k+1,k+2,\dotsi,2k\}$ and $x=\{i_1,i_2,\dotsi,i_k\},y=\{i_{k+1},i_{k+2},\dotsi,i_{2k}\}$ such that $i_{2k}=n$ and $i_p \ne i_q$ for all $p \ne q$. Since $2k=n \ge 6$, it is easy to find two vertices $a,b$ that are different from $v$ and $y$ in $Y$. Both $a$ and $b$ have $k(\ge 3)$ adjacent vertices in $X$. Thus, let $a'$, $b'(\ne u,x)$ be vertices adjacent to $a$ and $b$ in $X$, respectively. Thus, there are two paths of $P2C(u,a';x,b')$ in $X$ and two paths of $P2C(a,v;b,y)$ in $Y$ by induction. Concatenating them by edges $aa'$ and $bb'$, yields two paths of $P2C(u,v;x,y)$ in $J(n,k)$. This completes the proof.$\qed $

  \begin{lemma}
  ~Let $A= \{a_1,a_2,\dotsi,a_m\}$ and let $s$ and $t$ be any two distinct vertices in $QJ(n,A)$.

  (i) If $1 \le a_1<a_2<\dotsi<a_m < n-1$ and $n\ge 4$, then any two vertices in $J(n,a_i)$, $1 \le i< m$, each of them has a distinct neighbor in $J(n,a_{i+1})-\{s,t\}$.

  (ii) If $1 < a_1<a_2<\dotsi<a_m \le n-1$ and $n\ge 4$, then any two vertices in $J(n,a_i)$, $1 < i \le m$, each of them has a distinct neighbor in $J(n,a_{i-1})-\{s,t\}$.

  In particular, if $A=\{1, n-1\}$, then any two vertices in $J(n,1)$, each of them has a distinct neighbor in $J(n,n-1)-\{s,t\}$, and any two vertices in $J(n,n-1)$, each of them has a distinct neighbor in $J(n,1)-\{s,t\}$.
  \end{lemma}
  \noindent Proof. (i) We may assume that $p=a_i<a_{i+1}=q$, and $u=\{i_1,\dotsi,i_p\},v=\{j_1,\dotsi,j_p\}$ in $J(n,a_i)$. Thus, $u$ or $v$ has $\binom{n-p}{q-p}$ neighbors in $J(n,a_{i+1})$. Since $a_m <n-1$ and $i<m$, $a_i \le n-3$. It follows that $\binom{n-p}{q-p} \ge \binom{3}{1}=3$, meaning that $u$ or $v$ has at least three neighbors in $J(n,a_{i+1})$.
 Thus, there is at least one neighbor of $u$ (resp. $v$) in $J(n,a_{i+1})-\{s,t\}$.
 Let $u'$ and $v'$ be neighbors of $u$ and $v$ in $J(n,a_{i+1})-\{x,y\}$, respectively. If $u'\neq v'$, we are done.
 So we assume that $u'=v'$. So $u'$ is a $q$-subset that contains $i_1,\dotsi,i_p,j_1,\dotsi,j_p$. Since $q<n-1$, there are at least two elements $k_1, k_2\in [n]$ but $k_1,k_2\notin u'$. By replacing $i_1$ (resp. $j_1$) with $k_i(i=1,2)$ in $u'$, we obtain two vertices, say $a_1,a_2$ (resp. $b_1,b_2$). Thus, $a_1,a_2$ (resp. $b_1,b_2$) are neighbors of $v$ (resp. $u$). And these four vertices are different from each other. So, it is easy to find two distinct neighbors of $v$ and $u$ in $\{a_1,a_2,b_1,b_2,u'\}-\{s,t\}$, respectively.

(ii) As $J(n,k)$ is isomorphic to $J(n,n-k)$, this statement is equivalent to the following: any two vertices in $J(n,n-a_i)$, each of them has a distinct neighbor in $J(n,n-a_{i-1})-\{s,t\}$ where $s,t$ are two distinct vertices in $J(n,n-a_{i-1})$, since $J(n,k)$ is isomorphic to $J(n,n-k)$. Let $b_{m+1-i}=n-a_i$ for $i=1,2,\cdots , m$. It follows that $1\le b_1<b_2<\cdots <b_m<n-1$. Thus, we obtain that any two vertices in $J(n,b_j)$, $1\le j < m$, each of them has a distinct neighbor in $J(n,b_{j+1})-\{s,t\}$ by the proof of (i). Letting $j=m+1-i$, we can obtain that $b_j=n-a_i$ and $b_{j+1}=n-a_{i-1}$, meaning that the statement is true.

  In particular, if $A=\{1,n-1\}$, we may assume that $u=\{i\}$ and $v=\{j\}$ in $J(n,1)$. Clearly, $u$ or $v$ has $n-1$ neighbors in $J(n,n-1)$ as there is exactly
one $(n-1)$-subset that does not contain $i$ or $j$. It follows that $u$ and $v$ have $n-2$ neighbors in common. Hence, $u$ and $v$, each of them has a distinct neighbor in $J(n,n-1)-\{s,t\}$. In addition, any two vertices in $J(n,n-1)$, each of them has a distinct neighbor in $J(n,1)-\{s,t\}$ as $J(n,1)$ is isomorphic to $J(n,n-1)$. This completes the proof.$\qed$

 \begin{lemma}
  ~Let $A= \{a_1,a_2,\dotsi,a_m\}$ with $1 \le a_1<a_2<\dotsi<a_m \le n-1$ and $n\ge 4$, and let $s$ be an arbitrary vertex in $QJ(n,A)$. Then

  (i) any vertex in $J(n,a_i)$, $1 \le i< m$, has a neighbor in $J(n,a_{i+1})-\{s\}$;

  (ii) any vertex in $J(n,a_i)$, $1 < i\le m$, has a neighbor in $J(n,a_{i-1})-\{s\}$.
  \end{lemma}
 \noindent Proof. We may assume that $p=a_i<a_{i+1}=q$, and let $u$ be a vertex in $J(n,a_i)$. Thus, $u$ has $\binom{n-p}{q-p}$ neighbors in $J(n,a_{i+1})$. Since $a_m \le n-1$ and $i<m$, $a_i \le n-2$. It follows that $\binom{n-p}{q-p} \ge \binom{2}{1}=2$, meaning that $u$ has at least two neighbors in $J(n,a_{i+1})$. Thus, any vertex in $J(n,a_i)$ has a neighbor in $J(n,a_{i+1})-\{s\}$. Similarly, the statement (ii) also holds.$\qed$

 \begin{lemma}
     ~Let $A=\{a_1,a_2,\cdots,a_m\}$, $1\le a_1 < \cdots < a_m \le n-1$. If $QJ(n,A)$ is paired 2-coverable, then $QJ(n,A\cup \{n\})$ is paired 2-coverable.
 \end{lemma}
 \noindent Proof. Let $u,v,x$ and $y$ be four endpoints in $QJ(n,A\cup \{n\})$. Since $J(n,n)$ has only one vertex, i.e. $c=\{1,2,\cdots,n\}$, $c$ is adjacent to all vertices in $J(n,a_m)$.\\
 \textbf{Case 1}. $c \neq u,v,x,y$, meaning $u,v,x,y \in V(QJ(n,A))$ . Let $ab \in E(J(n,a_m))$ be an edge on one of the two paths of $P2C(u,v;x,y)$ in $QJ(n,A)$. Deleting $ab$ and adding edges $ac,bc$, yields the two paths of $P2C(u,v;x,y)$ in $QJ(n,A\cup \{n\})$.\\
 \textbf{Case 2}. $c$ is one of $u,v,x$ and $y$. We may assume that $c=u$. Choose a vertex $c'\neq v,x,y$ in $J(n,a_m)$. By our assumption, there are two paths $P$ and $Q$ of $P2C(c',v;x,y)$ in $QJ(n,A)$. Adding $cc'$ yields two paths of $P2C(u,v;x,y)$ in $QJ(n,A\cup \{n\})$. This completes the proof.\qed

 \begin{lemma}
     ~Let $A=\{a_1,a_2,\cdots,a_m\}$, $1\le a_1 < \cdots < a_m \le n-1$. Then $QJ(n,A)$ is paired 2-coverable.
 \end{lemma}
 \noindent Proof. As $QJ(n, \{a_p, a_{p+1},\cdots , a_{p+q}\})$, $1\le p \le p+q \le m$, is an induced subgraph of $QJ(n,A)$, we call two paths of $P2C(u,v;x,y)$ of $QJ(n, \{a_p, a_{p+1},\cdots , a_{p+q}\})$ {\em local $P2C$ paths} of $Q(n,A)$. Now we give a method to expand local $P2C$ paths to $P2C$ paths of $QJ(n,A)$. For simplicity, we call the method $EP2C$. If $p>1$, it is easy to find an edge $ab \in E(J(n,a_p))$ on one of two paths of $P2C(u,v;x,y)$ in $QJ(n,\{a_p,a_{p+1},\cdots,a_{p+q}\})$. Let $a'$ and $b'$ be neighbors of $a$ and $b$ in $J(n,a_{p-1})$, respectively. By Lemma 2, there is a Hamilton path $P_1$ in $QJ(n,\{a_1,\cdots ,a_{p-1}\})$ with endpoints $a'$ and $b'$. If $p+q<m$, it is easy to find an edge $cd \in E(J(n,p+q))$ on one of two paths of $P2C(u,v;x,y)$ in $QJ(n,\{a_p,a_{p+1},\cdots,a_{p+q}\})$, such that neither $a$ nor $b$ is on $cd$. Let $c'$ and $d'$ be neighbors of $c$ and $d$ in $J(n,a_{p+q+1})$, respectively. By Lemma 2, there is a Hamilton path $P_2$ in $QJ(n,\{a_{p+q+1},\cdots ,a_m\})$ with endpoints $c'$ and $d'$. Deleting edges $ab,cd$, adding $aa',bb',cc',dd'$ and concatenating $P_1,P_2$ and local $P2C$ paths, yields $P2C$ paths of $QJ(n,A)$. Thus, we can prove this lemma by extending local $P2C$ paths to $P2C$ paths of $QJ(n,A)$ by using $EP2C$.

 Next we distinguish following cases according to which levels the four endpoints locate.\\
 \textbf{Case 1}. $u,v,x$ and $y$ are contained in exactly one level. There are $P2C$ paths in this level by Theorem 4, which are local $P2C$ paths of $QJ(n,A)$. Thus, they can be extended to $P2C$ paths of $QJ(n,A)$ by $EP2C$.

\begin{figure}[h]
	\centering
	\includegraphics[width=0.7\textwidth]{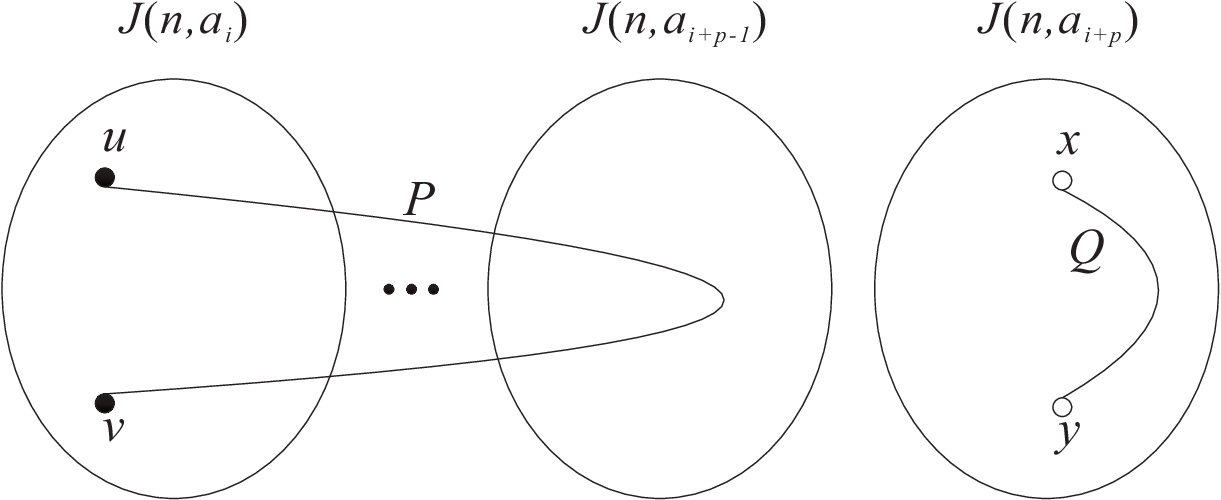}
     \caption{Local $P2C$ paths of Case 2.1. }
     \label{case2.1}
\end{figure}

  \noindent  \textbf{Case 2}. $u,v,x$ and $y$ are contained in two levels. We may denote the two levels by $i$ and $i+p$, where $p>0$.\\
\textbf{Case 2.1}. Level $i$ contains $u,v$ and level $i+p$ contains $x,y$.
There is a Hamilton path $P$ in $QJ(n,\{a_i,a_{i+1},\cdots,a_{i+p-1}\})$ with endpoints $u,v$ by Lemma 2. In addition, there is a Hamilton path $Q$ in $J(n,a_{i+p})$ with endpoints $x,y$ by Lemma 1. Thus, $P$ and $Q$ are local $P2C$ paths (see Fig.~\ref{case2.1}), which can be extended to $P2C$ paths of $QJ(n,A)$.

\begin{figure}
	\centering
	\includegraphics[width=0.7\textwidth]{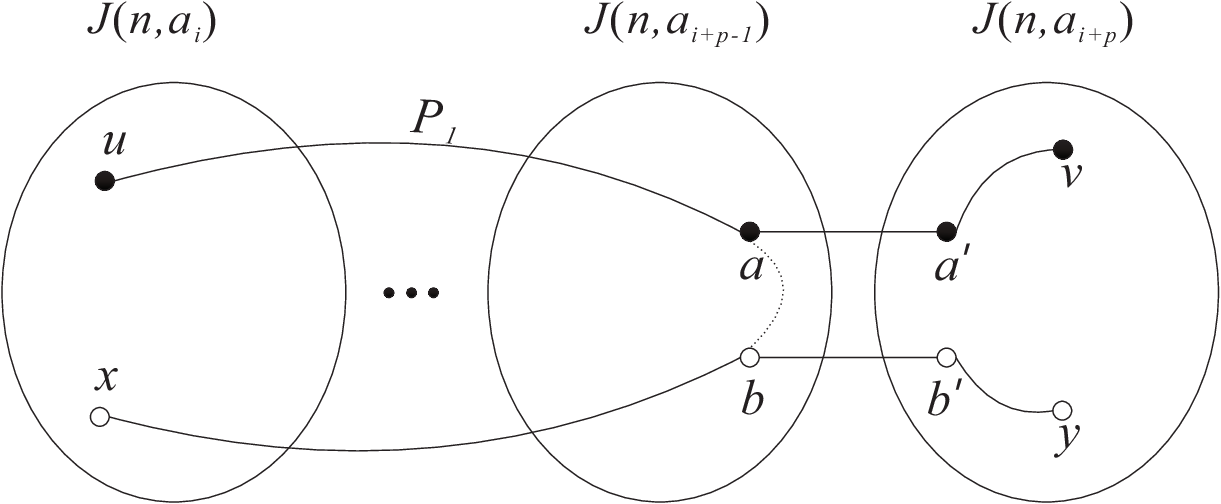}
     \caption{Local $P2C$ paths of Case 2.2 when $a_{a+p}<n-1$.}
     \label{case2.2.1}
\end{figure}

\begin{figure}
	\centering
	\includegraphics[width=0.7\textwidth]{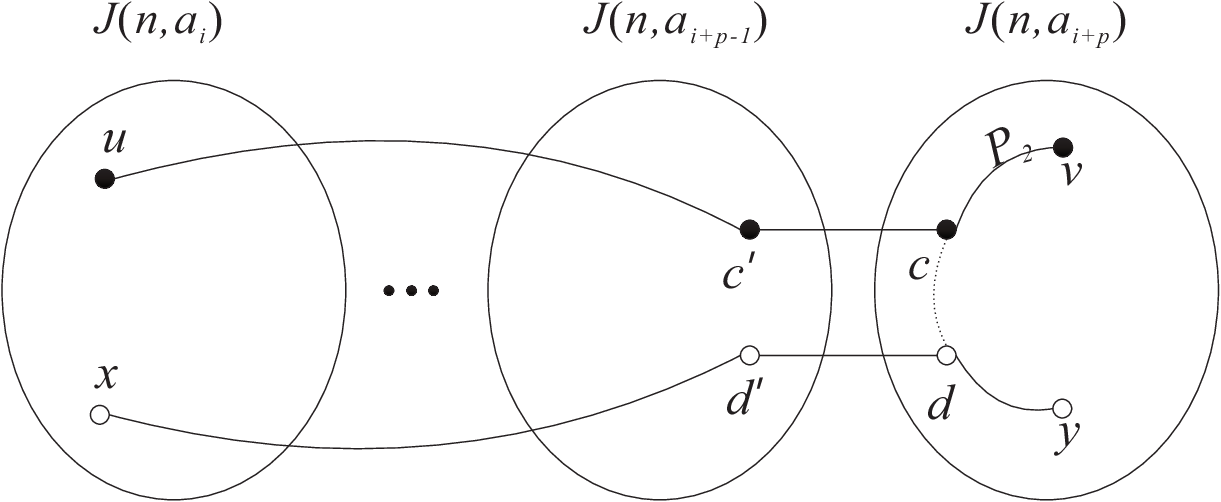}
     \caption{Local $P2C$ paths of Case 2.2 when $a_{a+p}=n-1$.}
     \label{case2.2.2}
\end{figure}

\noindent \textbf{Case 2.2}. Level $i$ contains $u,x$ and level $i+p$ contains $v,y$.
    If $a_{i+p}<n-1$, firstly, we can obtain a Hamilton path $P_1$ in $QJ(n,\{a_i,\cdots,a_{i+p-1}\})$ with endpoints $u$ and $x$ by Lemma 2.
 Let $ab \in E(J(n,a_{i+p-1}))$ be an edge on $P_1$. Then, by Lemma 5, we can choose $a'$ and $b'$ in $J(n,a_{i+p})-\{v,y\}$ as neighbors of $a$ and $b$, respectively.
 By Theorem 4, there are two paths of $P2C(a',v;b',y)$ in $J(n,a_{i+p})$.
 Deleting $ab$, adding edges $aa',bb'$ and concatenating these paths, yields local $P2C$ paths (see Fig.~\ref{case2.2.1}).
 If $a_{i+p}=n-1$, we can find a Hamilton path $P_2$ in $J(n,a_{i+p})$ with endpoints $v$ and $y$ by Lemma 1. Let $cd$ be an edge on $P_2$. Then, by Lemma 5, we can choose $c'$ and $d'$ in $J(n,a_{i+p-1})-\{u,x\}$ as neighbors of $c$ and $d$, respectively. Since $a_{i+p-1}<n-1$, there are two paths of $P2C(c',u;d',x)$ in $QJ(n,\{a_{i},\cdots,a_{i+p-1}\})$.
 Deleting $cd$, adding edges $cc',dd'$ and concatenating these paths, yields local $P2C$ paths (see Fig.~\ref{case2.2.2}). Then, we can obtain the $P2C$ paths of $QJ(n,A)$ by using $EP2C$.


\begin{figure}
	\centering
	\includegraphics[width=0.7\textwidth]{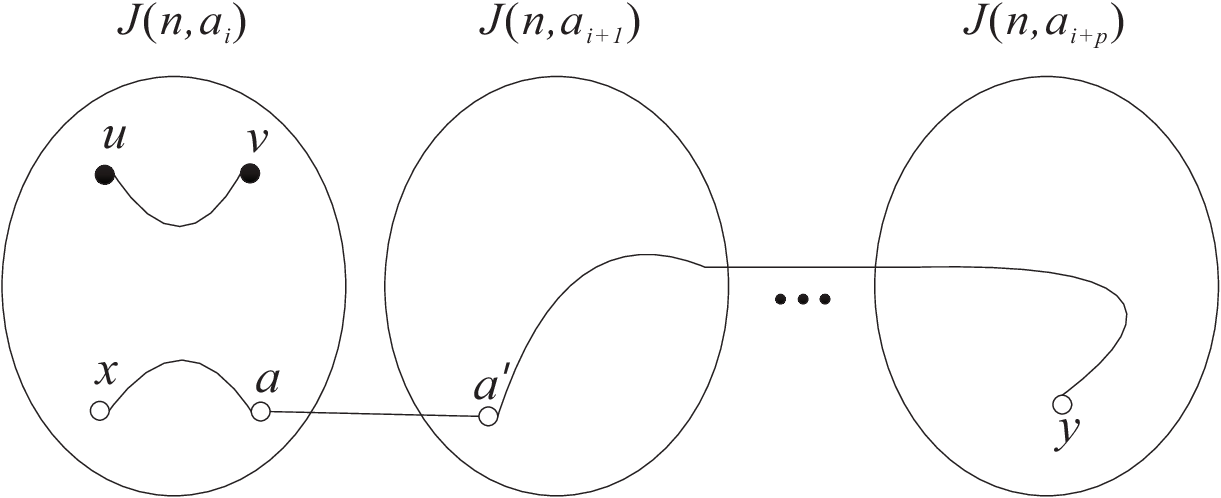}
     \caption{Local $P2C$ paths of Case 2.3. }
     \label{case2.3}
\end{figure}

\noindent \textbf{Case 2.3}. Three endpoints are contained in level $i$ (or $i+p$) and the other endpoint is contained in level $i+p$ (or $i$). Without loss of generality, we may assume that $u,v,x$ are contained in level $i$ and $y$ is contained in level $i+p$.
    It is clear that any vertex $a$ in $J(n,a_i)-\{u,v,x\}$ has a neighbor $a'$ in $J(n,a_{i+1})-\{y\}$ by Lemma 6. Then, there are two paths of $P2C(u,v;x,a)$ in $J(n,a_i)$ by Theorem 4, and a Hamilton path with endpoints $a',y$ in $QJ(n,\{a_{i+1},\cdots,a_{i+p}\})$ by Lemma 2. Adding $aa'$ and concatenating these paths, yields local $P2C$ paths (see Fig.~\ref{case2.3}). Similarly, we can obtain $P2C$ paths of $QJ(n,A)$.

\begin{figure}
	\centering
	\includegraphics[width=1\textwidth]{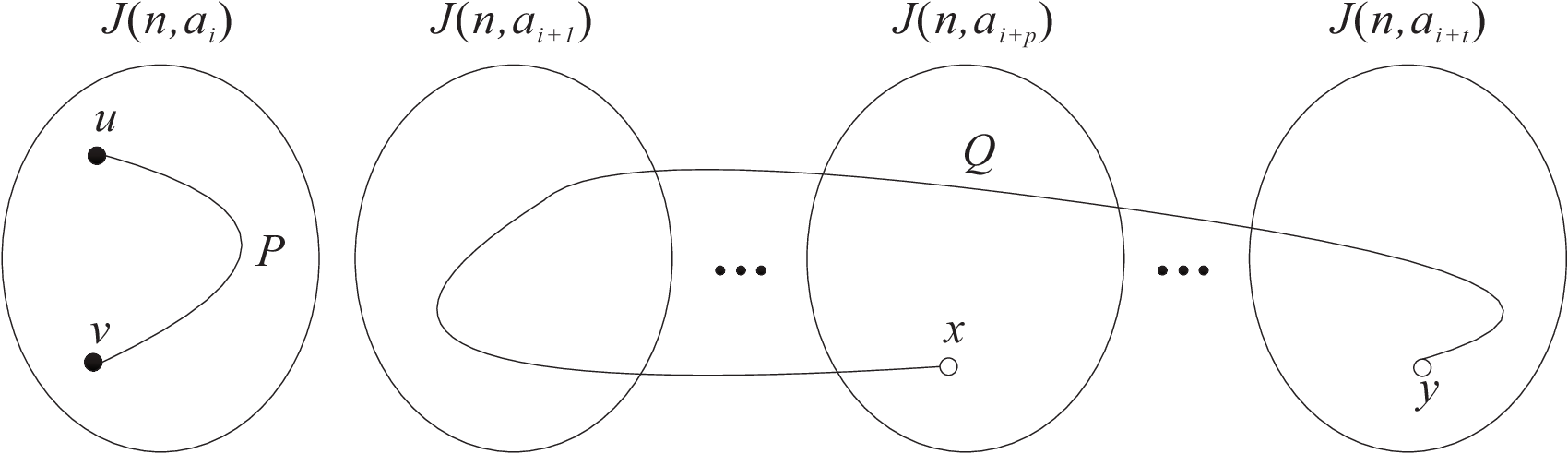}
     \caption{Local $P2C$ paths of Case 3.1.1. }
     \label{case3.1.1}
\end{figure}

 \noindent   \textbf{Case 3}. $u,v,x$ and $y$ are contained in exactly three levels. We may denote the three levels by $i,i+p$ and $i+t$, respectively, where $0<p<t$. It is clear that the case of level $i$ containing exactly two endpoints is equivalent to that of level $i+t$ containing exactly two endpoints. Thus, we further distinguish the following cases. \\
\textbf{Case 3.1}. Level $i$ contains two endpoints, and levels $i+p$ and $i+t$ contains exactly one endpoint, respectively.\\
\textbf{Case 3.1.1}. Level $i$ contains $u,v$, and level $i+p$ and level $i+t$ contain $x,y$, respectively. There are a Hamilton path $P$ with endpoints $u,v$ in $J(n,a_i)$, and a Hamilton path $Q$ with endpoints $x,y$ in $QJ(n,\{a_{i+1},\cdots,a_{i+t}\})$ by Lemma 2. Thus, $P$ and $Q$ are local $P2C$ paths of $QJ(n,A)$ (see Fig.~\ref{case3.1.1}), which can be extended to $P2C$ paths of $QJ(n,A)$ by $EP2C$.

\begin{figure}
	\centering
	\includegraphics[width=0.9\textwidth]{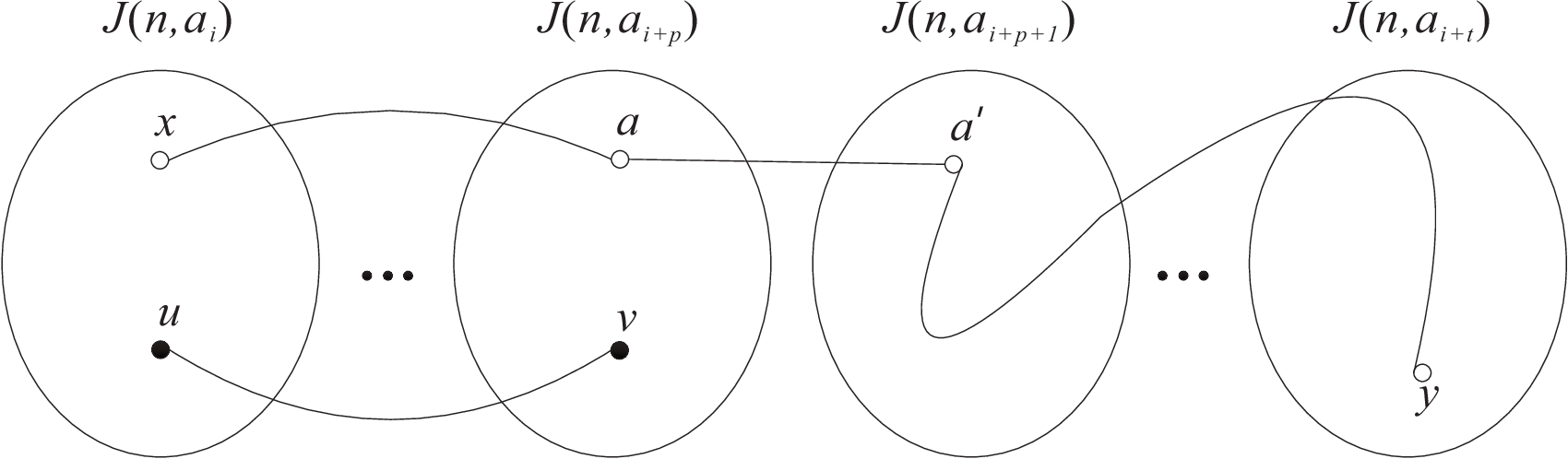}
     \caption{Local $P2C$ paths of Case 3.1.2. }
     \label{case3.1.2}
\end{figure}

\noindent \textbf{Case 3.1.2}. Level $i$ contains $u,x$, and level $i+p$ and $i+t$ contain $v$ and $y$, respectively.
Any vertex $a$ in $J(n,a_{i+p})-\{v\}$ has a neighbor $a'$ in $J(n,a_{i+p+1})-\{y\}$ by Lemma 6. Thus, there are two paths of $P2C(u,v;x,a)$ in $QJ(n,\{a_i,\cdots,a_{i+p}\})$ by the proof of Case 2.2 and a Hamilton path in $QJ(n,\{a_{i+p+1},\cdots,a_{i+t}\})$ with endpoints $a',v_2$ by Lemma 2. Adding $aa'$ and concatenating these paths, yields local $P2C$ paths of $QJ(n,A)$ (see Fig.~\ref{case3.1.2}), which can be extended to $P2C$ paths of $QJ(n,A)$ by $EP2C$.

\begin{figure}
	\centering
	\includegraphics[width=0.9\textwidth]{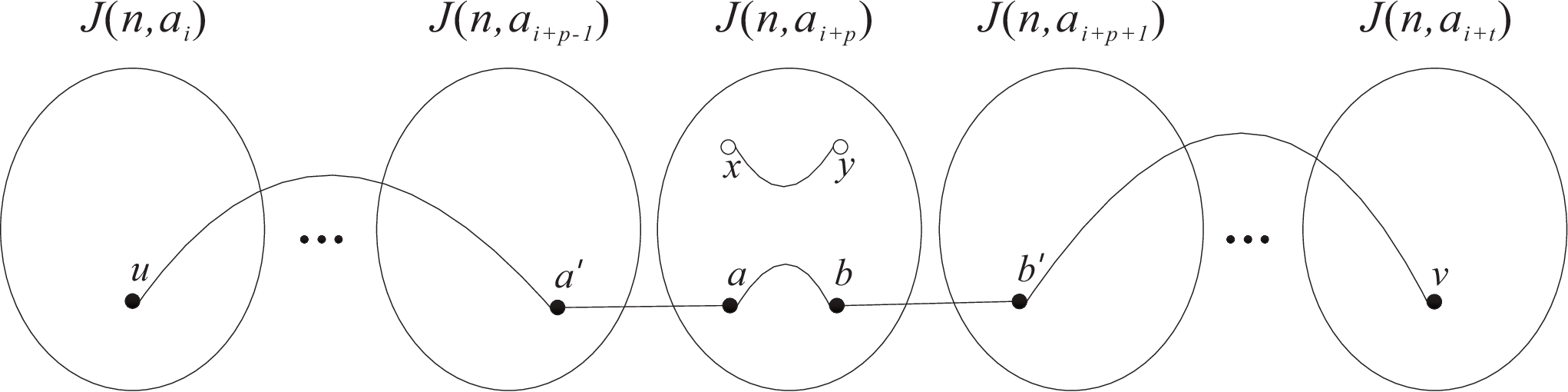}
     \caption{Local $P2C$ paths of Case 3.2.1. }
     \label{case3.2.1}
\end{figure}

\noindent \textbf{Case 3.2}. Level $i+p$ contains two endpoints, and levels $i$ and $i+t$ contains exactly one endpoint, respectively.\\
 \textbf{Case 3.2.1}. Level $i+p$ contains $x,y$, and level $i$ and $i+t$ contain $u$ and $v$,respectively.
 Let $a$ and $b$ be two distinct vertices in $J(n,a_{i+p})-\{x,y\}$. By Lemma 6, there exist a neighbor $a'$ of $a$ in $J(n,a_{i+p-1})-\{u\}$ and a neighbor $b'$ of $b$ in $J(n,a_{i+p+1})-\{v\}$. Thus, there are two paths of $P2C(x,y;a,b)$ in $J(n,a_{i+p})$ by Theorem 4. By Lemma 2, there are a Hamilton path in $QJ(n,\{a_i,\cdots,a_{i+p-1}\})$ with endpoints $u,a'$ and a Hamilton path in $QJ(n,\{a_{i+p+1},\cdots,a_{i+t}\})$ with endpoints $v,b'$. Adding edges $aa',bb'$, concatenating these paths, yields local $P2C$ paths (see Fig.~\ref{case3.2.1}). Similarly, we can obtain $P2C$ paths of $QJ(n,A)$.

\begin{figure}
	\centering
	\includegraphics[width=1\textwidth]{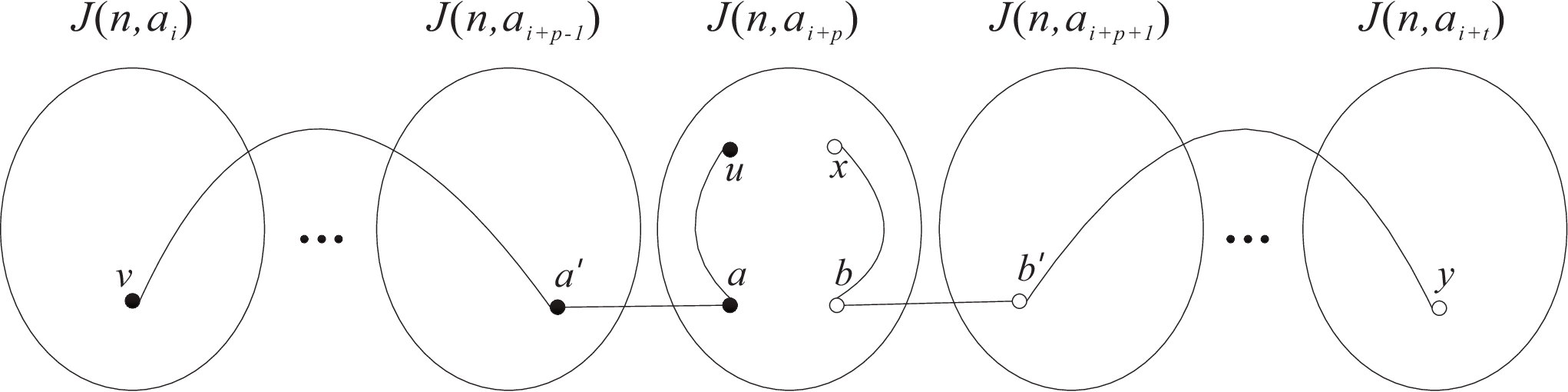}
     \caption{Local $P2C$ paths of Case 3.2.2. }
     \label{case3.2.2}
\end{figure}

\noindent \textbf{Case 3.2.2}. Level $i+p$ contains $u$ and $x$, and level $i$ and $i+t$ contain $v$ and $y$, respectively. The proof of this case is quite analogous to that of Case 3.2.1. Let $a$ and $b$ be two distinct vertices in $J(n,a_{i+p})-\{u,x\}$. Similarly, we choose a neighbor $a'$ of $a$ in $J(n,a_{i+p-1})-\{v\}$ and a neighbor $b'$ of $b$ in $J(n,a_{i+p+1})-\{y\}$. Thus, there are two paths of $P2C(u,a;x,b)$ in $J(n,a_{i+p})$ by Theorem 4. By Lemma 2, there are a Hamilton path in $QJ(n,\{a_i,\cdots,a_{i+p-1}\})$ with endpoints $v,a'$ and a Hamilton path in $QJ(n,\{a_{i+p+1},\cdots,a_{i+t}\})$ with endpoints $y,b'$. Adding edges $aa',bb'$ and concatenating these paths, yields local $P2C$ paths (see Fig.~\ref{case3.2.2}). Similarly, we can obtain $P2C$ paths of $QJ(n,A)$.

  \noindent  \textbf{Case 4}. $u,v,x$ and $y$ are contained in four different levels. We may denote the four levels by $i,i+p,i+s$ and $i+t$, where $0<p<s<t$.

\begin{figure}
	\centering
	\includegraphics[width=1\textwidth]{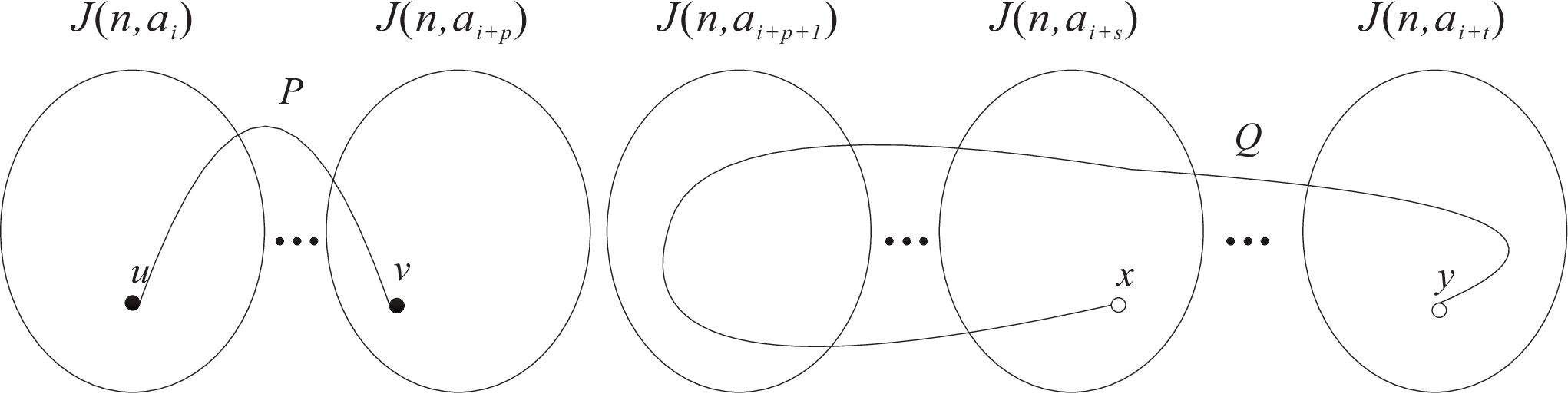}
     \caption{Local $P2C$ paths of Case 4.1. }
     \label{case4.1}
\end{figure}

\noindent \textbf{Case 4.1}. Levels $i,i+p,i+s$ and $i+t$ contain $u,v,x$ and $y$, respectively. There are a Hamilton path $P$ in $QJ(n,\{a_i,\cdots,a_{i+p}\})$ with endpoints $u,v$ and a Hamilton path $Q$ in $QJ(n,\{a_{i+p+1},a_{i+t}\})$ with endpoints $x,y$ by Lemma 2. Thus, $P$ and $Q$ are local $P2C$ paths of $QJ(n,A)$ (see Fig.~\ref{case4.1}), which can be extended to $P2C$ paths of $QJ(n,A)$ by $EP2C$.

\begin{figure}
	\centering
	\includegraphics[width=1\textwidth]{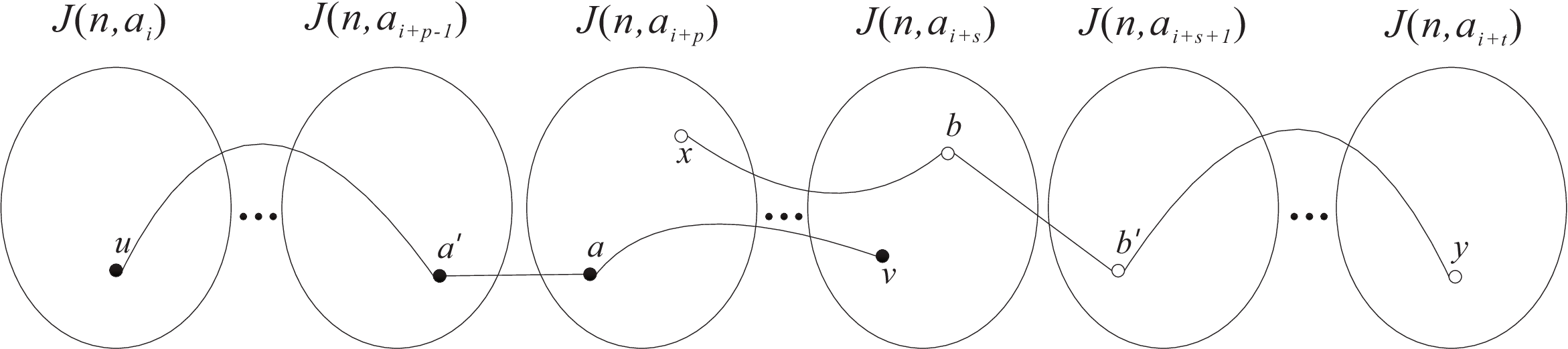}
     \caption{Local $P2C$ paths of Case 4.2. }
     \label{case4.2}
\end{figure}

\noindent \textbf{Case 4.2}. Levels $i,i+p,i+s,i+t$ contain $u,x,v$ and $y$, respectively. Let $a$ (resp. $b$) be a vertex in $J(n,a_{i+p})-\{x\}$ (resp. $J(n,a_{i+s})-\{v\}$). We can choose a neighbor $a'$ of $a$ in $J(n,a_{i+p-1})-\{u\}$ and a neighbor $b'$ of $b$ in $J(n,a_{i+s+1})-\{y\}$ by Lemma 6. Thus, there are two paths of $P2C(a,v;b,x)$ in $QJ(n,\{a_{i+p},\cdots,a_{i+s}\})$ by the proof of Case 2.2. In addition, by Lemma 2, there are a Hamilton path in $QJ(n,\{a_i,\cdots,a_{i+p-1}\})$ with endpoints $u,a'$ and a Hamilton path in $QJ(n,\{a_{i+s+1},\cdots,a_{i+t}\})$ with endpoints $y,b'$. Adding edges $aa',bb'$ and concatenating these paths, yields local $P2C$ paths (see Fig.~\ref{case4.2}). Similarly, we can obtain $P2C$ paths of $QJ(n,A)$ by using $EP2C$.
 \begin{figure}
	\centering
	\includegraphics[width=1\textwidth]{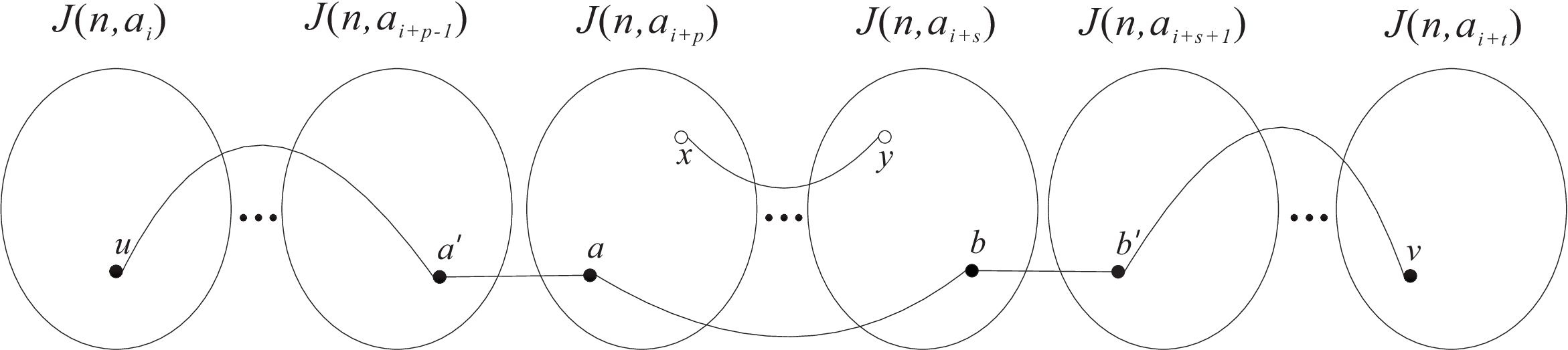}
     \caption{Local $P2C$ paths of Case 4.3. }
     \label{case4.3}
\end{figure}

\noindent \textbf{Case 4.3}. Levels $i,i+p,i+s,i+t$ contain $u,x,y$ and $v$, respectively. Let $a$ (resp. $b$) be a vertex in $J(n,a_{i+p})-\{x\}$ (resp. $J(n,a_{i+s})-\{y\}$). We can choose a neighbor $a'$ of $a$ in $J(n,a_{i+p-1})-\{u\}$ and a neighbor $b'$ of $b$ in $J(n,a_{i+s+1})-\{v\}$ by Lemma 6. Thus, there are two paths of $P2C(a,b;x,y)$ in $QJ(n,\{a_{i+p},\cdots,a_{i+s}\})$ by the proof Case 2.2. In addition, by Lemma 2, there are a Hamilton path in $QJ(n,\{a_i,\cdots,a_{i+p-1}\})$ with endpoints $u,a'$ and a Hamilton path in $QJ(n,\{a_{i+s+1},\cdots,a_{i+t}\})$ with endpoints $v,b'$. Adding edges $aa',bb'$ and concatenating these paths, yields local $P2C$ paths (see Fig.~\ref{case4.3}). Similarly, we can obtain $P2C$ paths of $QJ(n,A)$ by using $EP2C$.
    This completes the proof.$\qed$

\begin{theorem}
     ~Let $A$ be a nonempty subset of $[n]$. Then $QJ(n,A)$ is paired 2-coverable when $n\ge 4$ and $|V(QJ(n,A))|\ge 4$.
\end{theorem}
\noindent  Proof. If $n\notin A$, the theorem clearly holds by Lemma 8. If $n\in A$, let $A'=A-\{n\}$. Thus, $QJ(n,A')$ is paired 2-coverable by Lemma 8. It follows from Lemma 7 that $QJ(n,A)$ is paired 2-coverable. This completes the proof.$\qed$

\vskip 0.3 in

%

\noindent{\bf\normalsize Data Availability} Data sharing not applicable to this article as no datasets were generated or analysed during the current study.

\vskip 0.3 in

\noindent{\bf\large Declarations}

\vskip 0.3 in

\noindent{\bf\normalsize Competing Interests} The authors have not disclosed any competing interests.

\end{document}